\documentstyle{amsppt}
\voffset-10mm
\magnification1200
\pagewidth{130mm}
\pageheight{204mm}
\hfuzz=2.5pt\rightskip=0pt plus1pt
\binoppenalty=10000\relpenalty=10000\relax
\TagsOnRight
\loadbold
\nologo
\addto\tenpoint{\normalbaselineskip=1.2\normalbaselineskip\normalbaselines}
\addto\eightpoint{\normalbaselineskip=1.2\normalbaselineskip\normalbaselines}

\let\ge\geqslant
\let\wt\tilde
\let\[\lfloor
\let\]\rfloor
\redefine\d{\roman d}

\redefine\Re{\operatorname{Re}}

\define\ba{\boldkey a}
\define\bb{\boldkey b}
\define\bc{\boldkey c}

\define\fS{\frak S}
\define\fG{\frak G}
\define\fa{\frak a}
\define\fb{\frak b}
\define\fh{\frak h}
\define\pfrac#1#2{
\thickfracwithdelims..\thickness0{}{\lower2pt\rlap{\vrule height11.5pt}}
\kern-2pt
\frac{#1\,}{\,#2}
\kern-2.4pt
\thickfracwithdelims..\thickness0{\lower4.3pt\rlap{\vrule height11.5pt}}{}
}
\topmatter
\title
A few remarks on linear forms \\
involving Catalan's constant
\endtitle
\author
Wadim Zudilin \rm(Moscow)
\endauthor
\date
\hbox to100mm{\vbox{\hsize=100mm%
\centerline{E-print \tt math.NT/0210423}
\smallskip
\centerline{21 October 2002}
}}
\enddate
\address
\hbox to70mm{\vbox{\hsize=70mm%
\leftline{Moscow Lomonosov State University}
\leftline{Department of Mechanics and Mathematics}
\leftline{Vorobiovy Gory, GSP-2, Moscow 119992 RUSSIA}
\leftline{{\it URL\/}: \tt http://wain.mi.ras.ru/index.html}
}}
\endaddress
\email
{\tt wadim\@ips.ras.ru}
\endemail
\dedicatory
To N.\,M.~Korobov on the occasion of his 85th birthday
\enddedicatory
\abstract
In the joint work~\cite{RZ} of T.~Rivoal and the author,
a hypergeometric construction was proposed for studing
arithmetic properties of the values of Dirichlet's
beta function $\beta(s)$ at even positive integers.
The construction gives some bonuses \cite{RZ}, Section~9,
for Catalan's constant $G=\beta(2)$,
such as a second-order Ap\'ery-like recursion
and a permutation group in the sense of G.~Rhin and C.~Viola~\cite{RV}.
Here we prove expected integrality properties
of solutions to the above recursion as well as
suggest a simpler (also second-order and Ap\'ery-like) one for~$G$.
We `enlarge' the permutation group of~\cite{RZ}, Section~9,
by showing that the total $120$-permutation group of~\cite{RV}
for~$\zeta(2)$ can be applied in arithmetic study of Catalan's
constant. These considerations have computational meanings
and do not allow us to prove the (presumed) irrationality
of~$G$. Finally, we suggest a conjecture yielding
the irrationality property of numbers (e.g., of Catalan's constant)
from existence of suitable second-order difference equations
(recursions).
\endabstract
\endtopmatter
\footnote""{2000 {\it Mathematics Subject Classification}.\enspace
Primary 11J72; Secondary 11J70, 30B50, 33C60.}
\leftheadtext{W.~Zudilin}
\rightheadtext{Linear forms involving Catalan's constant}
\document

Recently, T.~Rivoal and the author~\cite{RZ} proved
several partial results on the irrationality of the numbers
$$
\beta(s)=\sum_{n=0}^\infty\frac{(-1)^n}{(2n+1)^s},
\qquad s=2,4,6,8,\dots\,.
$$
We did not succeed in proving the (expected) irrationality
of Catalan's constant $G=\beta(2)$. However, the general
analytic construction in~\cite{RZ} allows one to derive
a certain Ap\'ery-like second-order recursion
for Catalan's constant; this was done by semi-human
application of Zeilberger's creative telescoping in~\cite{Zu2}
and completely automatically, thanks to Ap\'ery's
`acc\'el\'eration de la convergence' approach, in~\cite{Ze}.

\medskip
I would like to thank S.~Fischler, T.~Rivoal, and J.~Sondow
for suggestions that allowed me to improve the text of the article.

\subhead
1. Hypergeometric series
\endsubhead
Recall the $\Bbb Q$-linear forms in $1$ and~$G$
constructed in \cite{RZ} and~\cite{Zu2}:
$$
\gather
{\align
r_n
&:=u_nG-v_n
\tag1
\\ &\phantom:
=\frac{n!}8\sum_{t=0}^\infty
(2t+n+1)\frac{\prod_{j=1}^n(t-j-1)\cdot\prod_{j=1}^n(t+j+n)}
{\bigl(\prod_{j=0}^n(t+j+\frac12)\bigr)^3}\,(-1)^t
\\ &\phantom:
=\frac{(-1)^nn!}8\,
\frac{\Gamma(3n+2)\,\Gamma(n+\frac12)^3\Gamma(n+1)}
{\Gamma(2n+\frac32)^3\Gamma(2n+1)}
\\ &\phantom:\qquad\times
{}_6\!F_5\biggl(\matrix\format&\,\r\\
3n+1, & \frac32n+\frac32, &  n+\frac12, &  n+\frac12, &  n+\frac12, &  n+1 \\
      & \frac32n+\frac12, & 2n+\frac32, & 2n+\frac32, & 2n+\frac32, & 2n+1
\endmatrix\biggm|-1\biggr),
\endalign}
\\ \vspace{1.5pt}
\lim_{n\to\infty}|r_n|^{1/n}
=\biggl|\frac{1-\sqrt5}2\biggr|^5,
\qquad
\lim_{n\to\infty}u_n^{1/n}
=\lim_{n\to\infty}v_n^{1/n}
=\biggl(\frac{1+\sqrt5}2\biggr)^5.
\endgather
$$
Note that Theorem~1 in~\cite{Zu2} states the following:
$$
2^{4n+3}D_nu_n\in\Bbb Z, \qquad 2^{4n+3}D_{2n-1}^3v_n\in\Bbb Z,
\tag2
$$
where $D_N$ denotes the least common multiple of
the numbers $1,2,\dots,N$, although the better inclusions
$$
2^{4n}u_n\in\Bbb Z, \qquad 2^{4n}D_{2n-1}^2v_n\in\Bbb Z
$$
hold for $n=1,2,\dots,1000$ by numerical verification
(see \cite{Zu2}, Section~4).
The aim of this section is to prove
(at least asymptotically, as $n\to\infty$, i.e., sufficient
for all practical purposes) this experimental observation.

\proclaim{Theorem 1}
For $n=0,1,2,\dots$, we have
$$
2^{4n+o(n)}u_n\in\Bbb Z,
\qquad 2^{4n+o(n)}D_{2n-1}^2v_n\in\Bbb Z.
\tag3
$$
\endproclaim

\remark{Remark}
As follows from the proof below,
the $o(n)$-term in the inclusions \thetag{3}
is of order $\log_2(2n)$.
\endremark

\demo{Proof}
We will require Whipple's transform
\cite{Ba}, Section~4.4, formula~(2),
$$
\align
&
{}_6\!F_5\biggl(\matrix\format&\,\c\\
a, & 1+\frac12a, &     b, &     c, &     d, &     e \\
   &   \frac12a, & 1+a-b, & 1+a-c, & 1+a-d, & 1+a-e
\endmatrix\biggm|-1\biggr)
\\ &\qquad
=\frac{\Gamma(1+a-d)\,\Gamma(1+a-e)}
{\Gamma(1+a)\,\Gamma(1+a-d-e)}
\cdot{}_3\!F_2\biggl(\gathered
1+a-b-c, \, d, \, e \\ 1+a-b, \, 1+a-c
\endgathered\biggm|1\biggr),
\tag4
\endalign
$$
provided that $\Re(1+a-d-e)>0$,
and Bailey's transform \cite{Ba}, Section~6.4, formula~(1),
$$
\align
&
{}_4\!F_3\biggl(\matrix\format&\,\c\\
a, \;\; & b, & c, & d \\ & k-b, & k-c, & k-d
\endmatrix\biggm|1\biggr)
\\ &\qquad
=\frac{\Gamma(k-b)\,\Gamma(k-c)\,\Gamma(k-d)}
{\Gamma(b)\,\Gamma(c)\,\Gamma(d)\,
\Gamma(k-b-c)\,\Gamma(k-b-d)\,\Gamma(k-c-d)}
\\ &\qquad\quad\times
\frac1{2\pi i}\int_{-i\infty}^{i\infty}
\frac{\aligned
\Gamma(b+t)\,\Gamma(c+t)\,\Gamma(d+t)\,\Gamma(k-a+2t)
\qquad \\ \vspace{-.4\baselineskip} \times
\Gamma(k-b-c-d-t)\,\Gamma(-t)
\endaligned}
{\Gamma(k-a+t)\,\Gamma(k+2t)}\,\d t,
\tag5
\endalign
$$
where the path of integration is parallel to the imaginary axis,
except that it is curved, if necessary, so that the decreasing
sequences of poles of the functions
$\Gamma(k-b-c-d-t)$ and $\Gamma(-t)$
lie to the left of the contour, while the increasing
sequences of poles of the functions
$\Gamma(b+t)$, $\Gamma(c+t)$, $\Gamma(d+t)$, and $\Gamma(k-a+2t)$
lie to the right.

Applying transform~\thetag{4} with
$a=3n+1$, $b=c=d=n+\frac12$, $e=n+1$
and then transform~\thetag{5} with
$a=2n+2$, $b=n+\frac12$, $c=d=n+1$, $k=3n+\frac52$
we obtain
$$
\align
r_n
&=\frac{(-1)^n}8\,\frac{(2n+1)!}{n!^2}
\frac1{2\pi i}\int_{-i\infty}^{i\infty}
\frac{\Gamma(n+1+t)^2\Gamma(n+\frac12+2t)\,\Gamma(-t)^2}
{\Gamma(3n+\frac52+2t)}\,\d t
\\
&=\frac{(-1)^n}8\,\frac{(2n+1)!}{n!^2}
\frac1{2\pi i}\int_{-i\infty}^{i\infty}
\frac{\Gamma(n+1+t)^2\Gamma(n+\frac12+2t)}
{\Gamma(1+t)^2\Gamma(3n+\frac52+2t)}
\biggl(\frac\pi{\sin\pi t}\biggr)^2\d t.
\tag6
\endalign
$$
Shifting $n+1+t\mapsto t$ and
considering the residues at the increasing sequence
of poles of the integrand in~\thetag{6}
(cf.~\cite{Ne}, Lemma~2) we arrive at the formula
$$
r_n
=\frac{(-1)^{n+1}}8\,\frac{(2n+1)!}{n!^2}
\sum_{\nu=1}^\infty
\frac{\d}{\d t}
\biggl(\frac{\Gamma(t)^2\Gamma(-n-\frac32+2t)}
{\Gamma(-n+t)^2\Gamma(n+\frac12+2t)}\biggr)\bigg|_{t=\nu}
=-\sum_{\nu=1}^\infty\frac{\d R_n(t)}{\d t}\bigg|_{t=\nu},
\tag7
$$
where
$$
\align
R_n(t)
&=\frac{(-1)^n}8\,\frac{(2n+1)!}{n!^2}
\frac{\bigl(\prod_{j=1}^n(t-j)\bigr)^2}
{\prod_{j=0}^{2n+1}(2t+j-n-\frac32)}
\\
&=\sum_{l=0}^n
\biggl(\frac{A_l}{t+l-\frac n2-\frac34}
+\frac{A_l'}{t+l-\frac n2-\frac14}\biggr).
\tag8
\endalign
$$
The coefficients $A_l$ and $A_l'$, $l=0,1,\dots,n$,
in the partial-fraction decomposition of the function~$R_n(t)$
can be easily determined by the standard procedure:
$$
\gather
\aligned
A_l
&=\frac{(-1)^n}{16}\,\frac{(2n+1)!}{(2l)!\,(2n-2l+1)!}
\cdot\biggl(\frac{(t-1)(t-2)\dotsb(t-n)}{n!}
\bigg|_{t=-l+\frac n2+\frac34}\biggr)^2,
\\
A_l'
&=\frac{(-1)^{n+1}}{16}\,\frac{(2n+1)!}{(2l+1)!\,(2n-2l)!}
\cdot\biggl(\frac{(t-1)(t-2)\dotsb(t-n)}{n!}
\bigg|_{t=-l+\frac n2+\frac14}\biggr)^2,
\endaligned
\\
l=0,1,\dots,n,
\endgather
$$
hence
$$
A_l\cdot2^{6n+4}\in\Bbb Z
\quad\text{and}\quad
A_l'\cdot2^{6n+4}\in\Bbb Z,
\qquad l=0,1,\dots,n,
\tag9
$$
by well-known properties of the integer-valued polynomials
$(t-1)(t-2)\dotsb(t-n)/n!$\,.
Using this decomposition we can continue formula~\thetag{7}
as follows:
$$
\align
r_n
&=16\sum_{l=0}^nA_l\cdot\sum_{\mu=0}^\infty\frac1{(4\mu+\epsilon)^2}
+16\sum_{l=0}^nA_l'\cdot\sum_{\mu=0}^\infty\frac1{(4\mu+\epsilon')^2}
\\ &\qquad
+16\sum_{l=0}^{m-1}A_l\sum_{\mu=1}^{m-l}\frac1{(4\mu-\epsilon)^2}
+16\sum_{l=0}^{m'-1}A_l'\sum_{\mu=1}^{m'-l}\frac1{(4\mu-\epsilon')^2}
\\ &\qquad
-16\sum_{l=m+1}^nA_l\sum_{\mu=0}^{l-m-1}\frac1{(4\mu+\epsilon)^2}
-16\sum_{l=m'+1}^nA_l'\sum_{\mu=0}^{l-m'-1}\frac1{(4\mu+\epsilon')^2},
\tag10
\endalign
$$
where $m=\[(n+1)/2\]$, $m'=\[n/2\]$;
$\epsilon=1$ for $n$~even and $\epsilon=3$ for $n$~odd;
$\epsilon'=4-\epsilon$;
$\[\,\cdot\,\]$ denotes the integer part of a number.
(For instance, if $n$~is even, we have
$$
\align
&
\sum_{\nu=1}^\infty\sum_{l=0}^n
\frac{A_l}{(\nu+l-\frac n2-\frac34)^2}
\\ &\qquad
=16\sum_{l=0}^{2m}A_l
\sum_{\nu=1}^\infty\frac1{(4(\nu+l-m-1)+1)^2}
=16\sum_{l=0}^{2m}A_l\sum_{\mu=l-m}^\infty\frac1{(4\mu+1)^2}
\allowdisplaybreak &\qquad
=16\sum_{l=0}^{m-1}A_l
\biggl(\sum_{\mu=l-m}^{-1}+\sum_{\mu=0}^\infty\biggr)\frac1{(4\mu+1)^2}
+16A_m\sum_{\mu=0}^\infty\frac1{(4\mu+1)^2}
\\ &\qquad\quad
+16\sum_{l=m+1}^{2m}A_l
\biggl(\sum_{\mu=0}^\infty-\sum_{\mu=0}^{l-m-1}\biggr)\frac1{(4\mu+1)^2}
\allowdisplaybreak &\qquad
=16\sum_{l=0}^{2m}A_l\cdot\sum_{\mu=0}^\infty\frac1{(4\mu+1)^2}
\\ &\qquad\quad
+16\sum_{l=0}^{m-1}A_l\sum_{\mu=1}^{m-l}\frac1{(4\mu-1)^2}
-16\sum_{l=m+1}^{2m}A_l\sum_{\mu=0}^{l-m-1}\frac1{(4\mu+1)^2}
\endalign
$$
and we proceed analogously in the three remaining cases.)

By~\thetag{8}, $R_n(t)=O(t^{-2})$ as $t\to\infty$,
hence $\sum_{l=0}^nA_l+\sum_{l=0}^nA_l'=0$.
Therefore, formula~\thetag{10} can be written
in the desired form $r_n=u_nG-v_n$, where
$$
\aligned
u_n
&=16(-1)^n\sum_{l=0}^nA_l
=16(-1)^{n+1}\sum_{l=0}^nA_l',
\\
v_n
&=-16\sum_{l=0}^{m-1}A_l\sum_{\mu=1}^{m-l}\frac1{(4\mu-\epsilon)^2}
-16\sum_{l=0}^{m'-1}A_l'\sum_{\mu=1}^{m'-l}\frac1{(4\mu-\epsilon')^2}
\\ &\qquad
+16\sum_{l=m+1}^nA_l\sum_{\mu=0}^{l-m-1}\frac1{(4\mu+\epsilon)^2}
+16\sum_{l=m'+1}^nA_l'\sum_{\mu=0}^{l-m'-1}\frac1{(4\mu+\epsilon')^2}.
\endaligned
\tag11
$$
(Again, Zeilberger's creative telescoping applied to the sequences
$u_n,v_n$ in~\thetag{11} yields the Ap\'ery-like
recursion from~\cite{Zu2} for the old sequences $u_n,v_n$
in~\thetag{1}; this fact implies the coincidence
of the two representations~\thetag{1} and~\thetag{11}
for the numbers $u_n$ and $v_n$ in the sequence $r_n=u_nG-v_n$.)
Formulae~\thetag{11} for $u_n,v_n$ and relations~\thetag{9}
imply $2^{6n}u_n,2^{6n}D_{2n-1}^2v_n\in\Bbb Z$. Finally,
using \thetag{2} and the fact that the order of~$2$ in~$D_N$
is $\[\log_2N\]$ we arrive at the desired inclusions~\thetag{3},
and the theorem is proved.
\enddemo

\subhead
2. A new Ap\'ery-like recursion for Catalan's constant
\endsubhead
The recursion in~\cite{RZ}, \cite{Zu2} allows one to do
fast computation of~$G$ with high accuracy.
Interpreting the solution to the recursion
in~\cite{RZ}, \cite{Zu2} as in~\thetag{7}
prompted us to modify slightly the parameters of the
above construction. Thus we take the sequence
$$
\gather
\wt r_n=\wt u_nG-\wt v_n
=-\sum_{\nu=1}^\infty\frac{\d\wt R_n(t)}{\d t}\bigg|_{t=\nu},
\tag12
\\
\wt R_n(t)
=\frac{(-1)^n}2\,\frac{(2n)!}{(n-1)!^2}
\frac{\prod_{j=1}^{n-1}(t-j)\cdot\prod_{j=1}^n(t-j)}
{\prod_{j=0}^{2n}(2t+j-n-\frac12)}
\endgather
$$
and apply Zeilberger's algorithm of creative telescoping
in order to prove the following result.

\proclaim{Theorem 2}
The numbers $\wt u_n$ and $\wt v_n$
satisfy the second-order recursion
$$
\align
&
(2n)^2(2n+1)^2(20n^2-20n+3)\wt u_{n+1}
-(3520n^6-2672n^4+196n^2-9)\wt u_n
\\ &\qquad
-(2n)^2(2n+1)(2n-3)(20n^2+20n+3)\wt u_{n-1}=0,
\qquad n=1,2,3,\dots,
\tag13
\endalign
$$
with the initial data
$\wt u_0=0$, $\wt u_1=6$,
and $\wt v_0=-1$, $\wt v_1=5$.
In addition, the limit relations
$$
\lim_{n\to\infty}|\wt u_nG-\wt v_n|^{1/n}
=\biggl|\frac{1-\sqrt5}2\biggr|^5,
\qquad
\lim_{n\to\infty}\wt u_n^{1/n}
=\lim_{n\to\infty}\wt v_n^{1/n}
=\biggl(\frac{1+\sqrt5}2\biggr)^5,
$$
hold and
$$
2^{4n+o(n)}\wt u_n\in\Bbb Z, \quad 2^{4n+o(n)}D_{2n-1}^2\wt v_n\in\Bbb Z
\qquad\text{for $n=0,1,2,\dots$}\,.
\tag14
$$
\endproclaim

(The proof of the inclusions \thetag{14} is a word-by-word repetition
of what was done in Section~1.)

The polynomials in~\thetag{13} are polynomials in~$2n$
with integer coefficients; the recursion~\thetag{13}
looks a little simpler than in~\cite{RZ},~\cite{Zu2}.
As in~\cite{Zu2}, Theorems~2 and~3, the constraint~\thetag{12}
also leads to the continued-fraction expansion
$$
\gather
6G=5+\pfrac{516}{q(2)}+\pfrac{p(3)}{q(4)}+\pfrac{p(5)}{q(6)}
+\dots+\pfrac{p(2n-1)}{q(2n)}+\dotsb,
\\
p(n)=(5n^2-20n+18)(n-2)(n-1)^2n^2(n+1)^2(n+2)(5n^2+20n+18),
\\
q(n)=55n^6-167n^4+49n^2-9,
\endgather
$$
and to the multiple Euler-type integral
$$
\gather
\wt u_nG-\wt v_n
=\frac{(-1)^{n-1}n}2\int_0^1\!\!\int_0^1
\frac{x^{n-3/2}(1-x)^ny^{n-1}(1-y)^{n-1/2}}
{(1-xy)^n}\,\d x\,\d y,
\\
n=1,2,3,\dots\,.
\endgather
$$

\subhead
3. A permutation group related to Catalan's constant
\endsubhead
Take the parameters $h_0,h_1,h_2,h_3,h_4$ satisfying
the conditions
$$
\gather
h_0,h_4\in\Bbb Z, \qquad h_1,h_2,h_3\in\Bbb Z+\tfrac12,
\tag15
\\
h_j>0 \quad\text{and}\quad 1+h_0-h_j-h_l>0
\qquad\text{for $j,l=1,2,3,4$}.
\tag16
\endgather
$$
As shown in~\cite{RZ}, Lemma~2, the quantity
$$
\align
&
\frac{\Gamma(1+h_0)\,\Gamma(h_3)\,\Gamma(h_4)\,
\Gamma(1+h_0-h_1-h_3)\,\Gamma(1+h_0-h_2-h_4)\,\Gamma(1+h_0-h_3-h_4)}
{\Gamma(1+h_0-h_1)\,\Gamma(1+h_0-h_2)\,
\Gamma(1+h_0-h_3)\,\Gamma(1+h_0-h_4)}
\\ &\qquad\times
{}_6\!F_5\biggl(\matrix\format&\,\c\\
h_0, & 1+\frac12h_0, & h_1, & h_2, & h_3, & h_4 \\
& \frac12h_0, & 1+h_0-h_1, & 1+h_0-h_2, & 1+h_0-h_3, & 1+h_0-h_4
\endmatrix\biggm|-1\biggr)
\tag17
\endalign
$$
belongs to the space $\Bbb QG+\Bbb Q$. (In~\cite{RZ},
a different ratio of gamma factors multiplies the
${}_6\!F_5$-series, but one ratio is a rational multiple
of the other.)

By means of the new parameters
$$
\gathered
a_1=1+h_0-h_1-h_2, \quad a_2=h_3, \quad a_3=h_4,
\\
b_2=1+h_0-h_1, \quad b_3=1+h_0-h_2
\endgathered
$$
and thanks to Whipple's transform~\thetag{4} we can represent
the quantity~\thetag{17} as follows:
$$
\align
&
\frac{\Gamma(a_2)\,\Gamma(a_3)\,\Gamma(b_2-a_2)\,\Gamma(b_3-a_3)}
{\Gamma(b_2)\,\Gamma(b_3)}
\cdot{}_3\!F_2\biggl(\matrix\format&\,\c\\
a_1, & a_2, & a_3 \\ & b_2, & b_3
\endmatrix\biggm|1\biggr)
\\ &\qquad
=\int_0^1\!\!\int_0^1
\frac{x^{a_2-1}(1-x)^{b_2-a_2-1}y^{a_3-1}(1-y)^{b_3-a_3-1}}
{(1-xy)^{a_0}}\,\d x\,\d y.
\tag18
\endalign
$$

Finally, take the third 10-element set~$\bc$:
$$
\gathered
c_{00}=(b_2+b_3)-(a_1+a_2+a_3)-1,
\\
c_{jl}=\cases
a_j-1 & \text{if $l=1$}, \\
b_l-a_j-1 & \text{if $l=2,3$},
\endcases
\qquad j,l=1,2,3
\endgathered
\tag19
$$
(hence all $c_{jl}>-1$ by~\thetag{16}),
in order to get that the double integral
$$
H(\bc)
=\int_0^1\!\!\int_0^1
\frac{x^{c_{21}}(1-x)^{c_{22}}y^{c_{31}}(1-y)^{c_{33}}}
{(1-xy)^{c_{11}+1}}\,\d x\,\d y
\tag20
$$
lies in $\Bbb QG+\Bbb Q$.
It will be useful to split the set~\thetag{19}
as $\bc=(\bc',\bc'')$, where
$$
\bc'=(c_{00},c_{21},c_{22},c_{33},c_{31})
\quad\text{and}\quad
\bc''=(c_{11},c_{23},c_{13},c_{12},c_{32})
$$
will be interpreted as cyclically ordered sets (i.e.,
$c_{00}$ follows $c_{31}$ in~$\bc'$ and
$c_{11}$ follows $c_{32}$ in~$\bc''$).
Obviously, each element in~$\bc''$ can be expressed in terms
of elements in~$\bc'$, and vice versa.
Using relations~\thetag{15} and summarizing what we said above
we obtain the following result.

\proclaim\nofrills{}
Suppose that
$$
c_{00},c_{21},c_{33}\in\Bbb Z+\tfrac12
\quad\text{and}\quad
c_{22},c_{31}\in\Bbb Z
\tag21
$$
for the elements in~$\bc'$
\rom(or, equivalently, $c_{13},c_{12},c_{32}\in\Bbb Z+\tfrac12$
and $c_{11},c_{23}\in\Bbb Z$ for the elements in~$\bc''$\rom)
and that all elements in~$\bc$ are $>-1$.
Then $H(\bc)\in\Bbb QG+\Bbb Q$.
\endproclaim

Digressing from the demi-integrality of the parameters~$\bc$,
let us note that the hypergeometric
${}_3\!F_2$-representation~\thetag{18} and
the equivalent ${}_6\!F_5$-representation~\thetag{17}
lead to the following group structure
(cf.~\cite{Wh} or \cite{Ba}, Sections~3.5--3.6).
Each permutation of the parameters $a_1,a_2,a_3$ in~\thetag{18}
or of the parameters $h_1,h_2,h_3,h_4$ in~\thetag{17}
gives a hypergeometric series of the same kind
(but with a different ratio of gamma factors before it).
For instance, the transposition $\fh=(h_1 \; h_4)$
rearranges the parameters $\ba$ and $\bb$ as follows:
$$
\fh\:\biggl(\matrix\format&\,\c\\
a_1, & a_2, & a_3 \\ & b_2, & b_3
\endmatrix\biggr)
\mapsto\biggl(\matrix\format&\,\c\\
b_3-a_3, & a_2, & b_3-a_1 \\ & b_2+b_3-a_1-a_3, & b_3
\endmatrix\biggr)
$$
and corresponds to Thomae's transformation \cite{Ba}, Section~3.2.
Hence the group~$\fG$ generated by all such permutations
appears naturally. An advantage of the superfluous $10$-element
set~$\bc$ is the fact that $\fG$
acts on the parameters~$\bc$ quite simply---by permutations.
As F.\,J.\,W.~Whipple has shown~\cite{Wh},
the group~$\fG$ is of order~$120$. A possible choice of
generators of~$\fG$ consists of the transpositions
$\fa_1=(a_1 \; a_3)$, $\fa_2=(a_2 \; a_3)$,
$\fb=(b_2 \; b_3)$, and the above-cited $\fh=(h_1 \; h_4)$
(see \cite{Zu1}, Section~6); the action of these permutations
on the set~$\bc$ reads as follows:
$$
\alignedat2
\fa_1&=(c_{11} \; c_{31})(c_{12} \; c_{32})(c_{13} \; c_{33}),
&\quad
\fa_2&=(c_{21} \; c_{31})(c_{22} \; c_{32})(c_{23} \; c_{33}),
\\
\fb&=(c_{12} \; c_{13})(c_{22} \; c_{23})(c_{32} \; c_{33}),
&\quad
\fh&=(c_{00} \; c_{22})(c_{11} \; c_{33})(c_{13} \; c_{31}).
\endalignedat
\tag22
$$

\proclaim{Theorem 3}
Let the quantity $H(\bc)$ be defined as the double integral
in~\thetag{20}, or as the ${}_3\!F_2$-series in~\thetag{18},
or as the ${}_6\!F_5$-series in~\thetag{17}.
Let $\fG\subset\fS_{10}$ be the $\bc$-permutation group
generated by~\thetag{22}.
Suppose that all elements in the set~$\bc$ are $>-1$.
\newline
Then
\roster
\item"(i)" the quantity
$$
\frac{H(\bc)}{\Pi(\bc)},
\qquad\text{where}\quad
\Pi(\bc)=\Gamma(c_{00})\,\Gamma(c_{21})\,
\Gamma(c_{22})\,\Gamma(c_{33})\,\Gamma(c_{31}),
\tag23
$$
is $\fG$-stable;
\item"(ii)" if the set~$\bc$ is $\fG$-equivalent to a set
satisfying condition~\thetag{21}, we have
$H(\bc)\in\allowmathbreak\Bbb QG+\Bbb Q$.
\endroster
\endproclaim

\demo{Proof}
(i) The $\fG$-stability of the quantity~\thetag{23}
has to be verified for the permutations in the list~\thetag{22};
this is routine using Whipple's transform
for verification of the $\fh$-stability.

(ii) In order to deduce the inclusion $H(\bc)\in\Bbb QG+\Bbb Q$
from the above claim~(i), it remains to show that
$\Pi(\sigma\bc)/\Pi(\bc)\in\Bbb Q$ for a set
$\bc$ satisfying~\thetag{21} and for all $\sigma\in\fG$
or, equivalently, for $\sigma\in\{\fa_1,\fa_2,\fb,\fh\}$
(by $\sigma\bc$ we mean the action of a permutation
$\sigma\in\fG$ on the set~$\bc$). This follows easily
from the fact that the gamma factors in
$$
\Pi(\bc), \quad \Pi(\fa_1\bc), \quad \Pi(\fa_2\bc), \quad
\Pi(\fb\bc), \quad \Pi(\fh\bc)
$$
have exactly three arguments belonging to $\Bbb Z+\frac12$ and two
arguments belonging to~$\Bbb Z$.
\enddemo

Another (very remarkable) description of the group~$\fG$
by means of the double integrals~\thetag{20} and their birational
transformations can be found in the work~\cite{RV}.

By~\cite{RV}, when all elements in~$\bc$ are
non-negative integers, one has
$H(\bc)\in\allowmathbreak\Bbb Q\zeta(2)+\Bbb Q$, where
$\zeta(2)=\pi^2/6$. Moreover, in this case,
$D_{m_1}D_{m_2}H(\bc)\in\Bbb Z\zeta(2)+\Bbb Z$,
where $m_1\ge m_2$ are the two successive maxima
of the set~$\bc$. This inclusion and the $\fG$-stability
of the quantity $H(\bc)/\Pi(\bc)$ make it possible to deduce
a nice irrationality measure for~$\zeta(2)$
(for details, see~\cite{RV}).

Theorems~1--3 allow us to expect a similar inclusion
$$
2^{2M+o(M)}D_{m_1}D_{m_2}H(\bc)\in\Bbb ZG+\Bbb Z
\tag24
$$
if the set $\bc$ is $\fG$-equivalent to a set
satisfying~\thetag{21}; here $M$~is the sum of two
integers in $\bc'=(c_{00},c_{21},c_{22},c_{33},c_{31})$
and $m_1\ge m_2$ are the two successive maxima of the set~$2\bc$.
Unfortunately, the inclusion~\thetag{24} is beyond the reach of
even the powerful group-structure approach to
proving irrationality results developed in~\cite{RV}
(see also \cite{Zu1}).

\subhead
4. Difference equations and irrationality
\endsubhead
Since
$$
\lim_{n\to\infty}D_{2n-1}^{1/n}=e^2
$$
by the prime number theorem, Theorem~1
(supplemented with equation~\thetag{1}) or Theorem~2
do not yield the irrationality of Catalan's constant.
What is the connection between irrationality
and Ap\'ery-like difference equations? We would like
to conclude this note by pointing out the following
expectation.

A sequence $\{x_n\}=\{x_n\}_{n=0}^\infty\subset\Bbb Q$
is said to satisfy the geometric condition\footnote{%
We should replace the standard term `$G$-condition'
by the phrase `geometric condition' since the capital letter~$G$
is reserved for Catalan's constant here.}
if the least common denominator of the numbers $x_0,x_1,\dots,x_n$
grows at most geometrically as $n\to\infty$.

Given a second-order recursion
$$
x_{n+1}+a(n)x_n+b(n)x_{n-1}=0,
\qquad
\lim_{n\to\infty}a(n)=a_0\in\Bbb Q,
\quad
\lim_{n\to\infty}b(n)=b_0\in\Bbb Q,
\tag25
$$
suppose that the characteristic polynomial $\lambda^2+a_0\lambda+b_0$
has roots $\lambda_1$ and $\lambda_2$
satisfying $0<|\lambda_1|<|\lambda_2|$.
Perron's theorem (see, e.g., \cite{Ge}, Chapter~V, Section~5)
then guarantees the existence of two linearly
independent solutions $\{x_n\}$ and $\{y_n\}$ such that
$$
\lim_{n\to\infty}\frac{x_{n+1}}{x_n}=\lambda_1,
\qquad
\lim_{n\to\infty}\frac{y_{n+1}}{y_n}=\lambda_2.
\tag26
$$

\proclaim{Conjecture}
In the above notation, suppose that both solutions
$\{x_n\}$ and $\{y_n\}$ of the recursion~\thetag{25}
are rational and satisfy the geometric condition.
Then $\lambda_1$ and $\lambda_2$ are rational numbers.
\endproclaim

This conjecture is trivially true
in the case of constant coefficients $a(n)=a_0$
and $b(n)=b_0$ of the recursion~\thetag{25};
we leave this observation as an exercise to the reader.

In order to show how the irrationality of $G$ follows from
the above conjecture, we have only to mention that,
if $G$~is rational, the solutions $\{\wt u_n\}$ and
$\{\wt r_n\}=\{\wt u_nG-\wt v_n\}$ to the recursion~\thetag{13}
are also rational numbers satisfying the geometric condition
and form Perron's basis,
while the roots $(11\pm5\sqrt5)/2=\bigl((1\pm\sqrt5)/2\bigr)^5$
of the characteristic polynomial are clearly irrational.

The geometric condition cannot be removed from hypothesis
of the conjecture\footnote{%
It is possible that the conjecture is true if we replace
the geometric condition hypothesis by the assumption
$a(n),b(n)\in\Bbb Q(n)$; however this new conjecture
would not cover several known cases
(for instance, the recursion corresponding
to Nesterenko's continued fraction for~$\zeta(3)$
in~\cite{Ne}, Theorem~2).}.
Indeed, taking $\lambda=(11+5\sqrt5)/2$
and $\lambda_1=-1/\lambda$, $\lambda_2=\lambda$, set
$$
x_n=\frac{(-1)^n}{\[\lambda^n\]}\in\Bbb Q,
\quad
y_n=\[\lambda^n\]\in\Bbb Z,
\qquad n=0,1,2,\dots\,.
\tag27
$$
Then $x_n\sim\lambda_1^n$ and $y_n\sim\lambda_2^n$ as $n\to\infty$,
hence relations~\thetag{26} hold. In addition, the sequences~\thetag{27}
satisfy the recursion~\thetag{25} with
$$
\aligned
b(n)&=-\frac{\[\lambda^{n-1}\]}{\[\lambda^{n+1}\]}
\cdot\frac{\[\lambda^n\]^2+\[\lambda^{n+1}\]^2}
{\[\lambda^{n-1}\]^2+\[\lambda^n\]^2},
\\
a(n)&=\frac{\[\lambda^n\]}{\[\lambda^{n-1}\]}\cdot b(n)
+\frac{\[\lambda^n\]}{\[\lambda^{n+1}\]},
\endaligned
\qquad n=0,1,2,\dots\,.
$$


\enddocument